\documentclass[a4paper,10pt,reqno]{amsart}

\usepackage{amsmath}
\usepackage{amsfonts,amssymb}
\usepackage{colordvi}
\usepackage{srcltx}
\usepackage{color}
\definecolor{Red}{rgb}{0.3,0.3,0.9}
\usepackage{graphicx,colordvi}

 \newcounter{zum}
 
\usepackage{amssymb, amsmath}

\usepackage{color}
\definecolor{Red}{rgb}{0.3,0.3,0.9}
\usepackage{multicol}
\DeclareMathOperator*{\esssup}{ess\,sup}

\begin{document}

\begin{center}
    { \large \bf On the Jackson constants for algebraic approximation of continuous functions}
\end{center}
\vskip .2cm
\begin{center}
    { A.\,G.\,Babenko, \ Yu.\,V.\,Kryakin}
\end{center}
\vskip .2cm
\begin{center}{\small \it Dedicated to Professor Igor A. Shevchuk on the occasion of his 70th birthday}
\end{center}



\ \

\vskip .2cm
\begin{quote}
\scriptsize {\bf Abstract}
{  We  establish { new estimates} for the constant $J_a(k,\alpha)$  in the Brudnyi--Jackson inequality for  
approximation  of $f \in C[-1,1]$ by algebraic polynomials:

$$
E_{n}^a (f) \le J_a(k, \alpha) \ \omega_k (f, \alpha \pi /n ),  \quad \alpha >0
$$
\vskip .2cm
\noindent
The main result of the paper implies the following {{inequalities}}

$$
1/2< J_a (2k, \alpha)  < 10, \quad  n \ge 2k(2k-1), \quad 
{ \alpha  \ge 2}
$$
}
\end{quote}

\section{Introduction}

\vskip .2cm
In this note we use the relatively  new approach (convolution series by C.Neumann \cite{n77} with the  Boman--Shapiro integral operators \cite{sh68, bosh71} ) for the constant problems in 
the following Brudnyi--Jackson theorem (see \cite{br68, dl93} ) for algebraic approximation of $f \in C[-1,1]$:
\begin{equation}\label{jb}
E_{n}^a (f) \le J_a(k, \alpha) \ \omega_k (f, \alpha \pi /n ),  \quad \alpha >0.
\end{equation}

The case of algebraic approximation  $E_{n-1}^a (f)$  of a  continuous  function  $f$   by algebraic po\-ly\-no\-mials of degree $ \le n-1$ is in some sense more
difficult then the case of trigonometric  approximation. 
Usually the reduction to the trigonometric approximation is used. There are  some  technical problems  in the case of the modulus of 
$\omega_k (f, \alpha \pi /n )$ of the order $k \ge 2.$
The problem of exact  constants in this case is a difficult one and we do not have  sharp results for $k > 1$.

We recall here  a result by Korneichuk \cite{kor63}, as the corollary of his remarkable theorem on constants in the case of concave modulus 
of continuity $\overline \omega $:
\begin{equation*}
E_{n-1}^a (f) \le  \frac 12 \ \overline \omega ( f, \pi/n),
\end{equation*}
\noindent
and the new Mironenko's result \cite{mir10} for the second modulus of continuity:

\begin{equation*}\label{mir}
E_{n-1}^a (f) \le  5 \   \omega_2 ( f, \alpha \pi/n), \quad \alpha=8^{-1/2}.
\end{equation*}

In the present paper we prove that { for $n>2k(2k-1)$} the constants in \eqref{jb} are bounded by an absolute { constant:}

\begin{equation*}\label{mr}
  {J_a} (2k, \alpha)  < 10, \quad  
{ \alpha \ge 2}.
\end{equation*}
{ It is clear that
$$
\omega_{2k+1}(f, \delta) \le 2 \ \omega_{2k}(f, \delta), \quad \delta >0,
$$
and therefore the main result of this paper states that  for $ \alpha \ge 2, \ n >k(k-1)$
we can write constant in \eqref{jb} that do not depend on $k$. }
\vskip .2cm

\section{Notation. Auxiliary facts. Main results}

In this paper  $i,j,k,l,m,n$  denote the natural numbers.
Let $A$ be  $\mathbb I = [-1,1]$ or $\mathbb R$.
We  denote by  {$W^k_*(A)$} the space of smooth functions $f^{(j)} \in C(A), \ j=0,1 \dots, k-1$, $f^{(k)}$ 
bounded a.e. on   $A$. We will also use the notation 
$$
\|f \| := \| f \|_A := \esssup_{x \in A} |f(x)|.
$$

We consider the approximation of real functions on  $\mathbb I =[-1,1]$ 
by algebraic polynomials $p_{n-1}(x) = \sum_{j=0}^{n-1} a_j x^j $ of degree at most  $n-1$. 
 We will denote by  $\mathcal P_{n-1}$ the space of such polynomials. 
The best approximation of $f \in C(\mathbb I)$ by  $p \in \mathcal P_{n-1}$ is defined by standard way
$$
E_{n-1}^a (f) := \inf_{p_{n-1}} \sup_{x \in \mathbb I} | f(x) - p_{n-1}(x)| = \inf_{p \in \mathcal P_{n-1} } \| f - p \|.
$$
Smoothness of function  $f \in C(A)$  is measured by modulus of {smoothness}.
Beside  the classical   $k$--th modulus of smoothness 
$$
\omega_k(f, \delta) := \sup_{x \in (1-kh/2)I, \ {0<h \le \delta}} |  \widehat\Delta_h^k f(x) |, \quad \widehat\Delta_h^k f(x):=\sum_{j=0}^k (-1)^{k-j} \binom{k}{j} f(x+jh-kh/2),
$$
we will use the special Boman--Shapiro modulus of continuity,
which measures the deviation of the function from the special linear combination of Steklov's means (see \cite{ st24, fks09, bks13}). 
\noindent
We will use the following convolution notation 
$$
(f*g)(x):= \int_{\mathbb{R}} f(t) \, g(x-t) \, dt,
$$
{ and the following notation for characteristic function and the convolution square of a characteristic function:
 $$
\chi_h (x) :=  \begin{cases}  \frac 1 h  &\quad x \in [-h/2,h/2],\\ 0, &\quad x \notin [-h/2,h/2].\end{cases}
$$
}
$$
\chi_h^2 (x) := (\chi_h * \chi_h )(x) = \begin{cases}  \frac 1 h {(1- |x|/h)}, &\quad x \in [-h,h],\\ 0, &\quad x \notin [-h,h].\end{cases}
$$
{ Define the special difference operator for a locally integrable function $f$  in the following way (see \cite{fks09, bks13})  }
$$
W_{2k}(f, x, \chi_h^2):=(-1)^k \frac 1{\binom {2k} k}\int_{\mathbb R} \widehat \Delta_t^{2k} f(x) \chi_h^2 (t) \, dt = (f - \Lambda_{2k}*f)(x),
$$
where
$$
\Lambda_{2k} (x) : = \Lambda_{2k,h}(x) = 2 { \sum_{j=1}^k }(-1)^{j+1} a_{j,k} \chi_{jh}^2 (x), \quad a_{j,k} := \binom{2k}{k+j}/\binom{2k}k.
$$
\noindent
{ It  was  proved in \cite{bks13}  that
\begin{equation*}\label{l2k}
\int_{\mathbb R} | \Lambda_{2k}(t) | \, dt \le  c_k -1,
\end{equation*}
with
$$
c_1 = 2, \ c_2 < 2.18, \ c_3< 2.26, \ c_4  < 2.31, \quad c_k < 3, \quad  k \ge 5,
$$
and therefore
\begin{equation}\label{bksb}
| W_{2k}(f, x, \chi_h^2)| \le c_k \, \sup_{t \in [x-kh, x+kh]} |f(t)|.
\end{equation}
}
\vskip .2cm
\noindent
We will use the standard notation for the Favard constants
$$\
\mathcal K_{k}:= \frac 4 {\pi} \, \sum_{j=-\infty}^\infty (4j +1)^{-k-1},
$$
$$
\mathcal K_0=1, \ \mathcal K_1= \frac \pi 2, \ \mathcal K_2 = \frac{\pi^2}{8}, \ \mathcal K_3 =\frac{\pi^3}{24},   \ \mathcal K_4 = \frac{5 \pi^4}{384}, \quad
\mathcal K_6 = \frac{61 \, \pi^6}{46080}, \quad \mathcal K_8 = \frac{277 \, \pi^8 }{2064384}.
$$

\vskip .2cm
We are now ready to state the main results of this paper.

\vskip .2cm
{\bf Theorem 1} {\it Suppose $f  \in C(\mathbb I)\,.$ Then {for $n \ge 2k(2k-1), \ k \ge 5$ } 
$$
E_{n-1}^a (f) \le  J_a(2k,\alpha) \, \omega_{2k} (f, \alpha \pi /n), 
$$
and 
$$
1/2 \le J_a(2k, \alpha) <  3 \, (2+e^{-2}) \,  \left(2 \ \sec(\pi/2\alpha) -1 - \frac{\mathcal K_2}{\alpha^2} \right),
\quad { 1<\alpha \le (2k-1)\pi^{-1}}.
$$
}
Note that in the trigonometric case we have (see \cite{fks09, bks13})
$$
{ {\binom{2k} k}^{-1} \le J(2k,\alpha) < \sec(\pi/2\alpha) \, {\binom{2k} k}^{-1} }, \quad \alpha > 1.
$$

\vskip .4cm
{\bf Theorem 2} {\it For small $k=1,2,3,4, \ \alpha > 0, $ we have the following estimates of the constants 
{
\begin{align*}
\frac 12 \le J_a(2, \alpha) &\le  { \frac 34 \left(1 + \frac{1}{4 \alpha^2} \right)}, \quad n \ge 2, \\
\frac 12 \le J_a(4, \alpha ) &\le 2.18 \cdot \left( 1 +  \mathcal K_2 \, \beta_2  \right) +  2 \, \mathcal K_4 \, \beta_2^2 ,  \quad n \ge  12, \\
\frac 12 \le J_a(6,\alpha) &\le 2.26 \cdot  ( 1 +  \mathcal{K}_2 \beta_3 + 2 \mathcal{K}_4 \beta_3^2 )  +  2 \mathcal{K}_6 \, \beta_3^3, \quad n \ge 30, \\
\frac 12 \le J_a(8,\alpha) &\le 2.31 \cdot  ( 1 +  \mathcal{K}_2 \beta_4 + 2 \mathcal{K}_4 \beta_4^2 + 2 \mathcal{K}_6 \beta_4^3 ) +  2 \mathcal{K}_8 \, \beta_4^4, \quad
n \ge 56,
\end{align*}
}
where
$$
\beta_{k} = \frac{4 \gamma_{k}}{\pi^2 \alpha^2},  \quad \gamma_2 = \frac{4}{3}, \quad \gamma_3 = \frac{68}{45}, \quad \gamma_4=
\frac{512}{315}.
$$
}
\vskip .2cm
Thus we may write the estimates  for 
 $\alpha =1,2$.  The constants  are increasing as $k$  increases for  $\alpha = 1$:
$$
J_a(2,1) \le  0.94, \quad  J_a(4,1) \le 4.38, \quad J_a(6,1) \le 6.71, \quad J_a(8,1) \le 8.9,
$$
but for 
$\alpha = 2$ 
$$
J_a(2,2) \le  0.8, \quad  J_a(4,2) \le 2.59, \quad J_a(6,2) \le 2.84, \quad J_a(8,2) \le 2.97,
$$
and for  $k \ge 5 $ the estimates of Theorem 1 give
$$
J_a(2k, 2) \le  9.74,   { \quad n  \ge 2k(2k-1)}.
$$
Note that if the Sendov conjecture $
w_k \le 1
$ is true  (see Theorem B below),
 then we achieve better inequality 
$$
J_a(2k,2) \le 9.74 \cdot (2+e^{-2})^{-1} = 4.57,
$$
which is near the results in the case of small $k$.

\vskip .2cm
{ The proofs of the main inequalities are based on the following fundamental facts. }
The first important fact is the algebraic variant \cite{sin81}  of the Bohr--Favard--Akhiesier--Krein inequality (see \cite{b35, f36, f37, ak37}):

\vskip .2cm
{\bf Theorem A} {\it  For $f \in  W_*^m (\mathbb I)$  
$$ 
E_{n-1}^a (f) \le  \frac {\mathcal K_m}{n^m} \, \| f^{(m)} \|, \quad m=1,2.
$$
$$
E_{n-1}^a (f) \le \mathcal{K}_m \  \frac{(n-m)!}{n!} \, \| f^{(m)} \|, \quad m \ge 3.
$$
}

\vskip .2cm
The second important fact is the modern variant of Whitney's theorem, with  good estimates of constants
(see \cite{w57,sen82, sen85, kry94, kry97,  kry02, gks02, zh02, zh04, ds08}):

\vskip .2cm
{\bf Theorem B} {\it { Suppose $f \in C[a,b], \ 0<a<b$. Then }
\begin{equation}\label{whb}
E_{k-1}^a (f)_\le \ w_k \ \omega_k \left(f, (b-a)/k \right),  
\end{equation}
with
\begin{equation*}
w_k \le 
\begin{cases}
 0.5, \quad &k=1,2,\\ 
 1, \quad &k =3, \dots 8,\\
 2, \quad &k \le 82 000, \\
 {  2 + \exp(-2)},  \quad &k \ge 82 000.
\end{cases}
\end{equation*}
}

{ Theorems A,B}  are the main technical tools for 
proving  Theorem 1 and Theorem 2.  Theorem 1 and Theorem 2 follow from  Proposition 1 and Proposition 2. 

\vskip .2cm

Proposition 1  concerns  the problem of  continuation of a function from $ \mathbb I$ to $ \mathbb R$. 
This continuation allows us to use the periodic--case approach.
\vskip .2cm
{\bf Proposition  1 }
{\it Let $f \in C(\mathbb I), \, 0<h< (2k)^{-1}. $
Then  there exists { a}  function  $g_f$ which is equal to  $f$ on  $\mathbb I$, continuous on $\mathbb{R}\setminus \mathbb I$ and such that
$$
W_{2k}(g_f,h):=\|  W_{2k}(g_f, \cdot, { \chi_h^2}) \|_{\mathbb R}  \le d_k \, \omega_{2k} (f, h),
$$
$$
 \|\widehat \Delta_h^{2k} g_f  \|_{\mathbb R}   \le d_k^*  \, \omega_{2k} (f, h),
$$
with
\begin{equation*}
d_k \le 
\begin{cases}
 1, \quad &k=1, \\
 c_k, \quad &k = 2,3, 4,\\
 6, \quad & 4 < k \le 41 000, \\
 3 (2 + \exp({-2})),  \quad &k > 41 000,
\end{cases}
\end{equation*}
where 
$$
c_2 < 2.18, \ c_3< 2.26, \ c_4  < 2.31,
$$
and 
\begin{equation*}
d_k^* \le 
\begin{cases}
 1.5, \quad &k =1, \\
 2^{2k}, \quad &k =2,3, 4,\\
 2^{2k+1}, \quad & 4 < k \le 41 000, \\
 {(2 +\exp (-2))  }\cdot 2^{2k},  \quad &k > 41 000.
\end{cases}
\end{equation*}

}

The key-idea is the same as in the periodic case (see  \cite{fks09,bks13}): we will use the truncated Neumann convolution series
{$$
g_f = \sum_{j=0}^{k-1} (g_f - g_f*\Lambda_{2k})*\Lambda_{2k}^j +  g_f*\Lambda_{2k}^{k}, \qquad \Lambda_{2k}^j:= \Lambda_{2k}^{j-1}*\Lambda_{2k},
$$
} 
with some  modification for the algebraic case.

\vskip .2cm
{\bf Proposition 2} {\it Let  $f \in C(\mathbb I)$, $g_f \equiv f$ on $\mathbb I$ and $g_f \in C( \mathbb R \setminus \mathbb I)$. { Then  for {$k \ge 2, \ n \ge 2k(2k-1)$}} 
\begin{align} \label{pr2}
E_{n-1}^a(f)
&\le\sum_{j=0}^{1} \frac{\mathcal{K}_{2j}}{n^{2j}} \left\|((g_f - g_f*\Lambda_{2k})*\Lambda_{2k}^j)^{(2j) } \right\|_{\mathbb R}
+ 2 \, \sum_{j=2}^{k-1}  \frac{\mathcal{K}_{2j}}{n^{2j}} \left\| ((g_f - g_f*\Lambda_{2k})*\Lambda_{2k}^j)^{(2j)} \right\|_{\mathbb R} \notag \\
& + 2 \, \frac{\mathcal{K}_{2k}}{n^{2k}} \, \left\| (g_f*\Lambda_{2k}^k)^{(2k)} \right\|_{\mathbb R}.
\end{align}
\vskip .2cm
}

{ For deducing the main theorems  from Proposition~2 we present here a new variant of the known estimates  (see \cite{fks09, bks13})}. 
\vskip .2cm

Put

\begin{equation*}\label{gamma}
 \gamma_k := 2 \sum_{j=1}^{(k+1)/2}  \frac{a_{2j-1,k}}{(2j-1)^2} < 2\sum_{j=1}^{\infty} \frac 1{(2j-1)^2} = \frac{\pi^2}4 , \quad  { k \ge 2.}
\end{equation*}

{\bf Lemma 1 } {\it { For  $ k \ge 2, \ j=0,\dots, k-1$ }
\begin{align*}
 { \| ((g_f - g_f*\Lambda_{2k})*\Lambda_{2k}^j)^{(2j) } \|_{\mathbb R}}  &\le \left( 4 \gamma_k \right)^j \, h^{-2j} \,  { \| g_f - g_f*\Lambda_{2k} \|_{\mathbb R}}, \\
 \left\| (g_f*\Lambda_{2k}^k)^{(2k)} \right\|_{\mathbb R} &\le { \gamma_k^k} \, h^{-2k} \| \widehat{\Delta}_h^{2k} g_f \|_{\mathbb{R}}. 
\end{align*}

}
Note that in  comparison with \cite{fks09,bks13} we add new inequality  here, which allows to estimate the last term in inequality \eqref{pr2} .

\vskip .2cm
To prove { Lemma 1} we represent 
$
\Lambda_{2k} (x)
$
{ as the special linear combination of }
\begin{equation*}\label{chi-s}
\chi_{h_j, h}^2 (x) :=\chi_h^2(x-h_j).
\end{equation*}

\vskip .2cm
{\bf Lemma  2} \ {{ \it For $k\ge 2, \ h > 0 $ we have}
\begin{equation*}\label{l21}
\Lambda_{2k}= 2 { \sum_{j=1}^k} (-1)^{j+1} a_{j,k} \chi_{jh}^2  = \sum_{j=-(k-1)}^{k-1} \alpha_{j,k} \, 
\chi_{jh, h}^2,  \quad \sum_{j=-k+1}^{k-1} | \alpha_{j,k} | = \gamma_k.
\end{equation*}
}
\vskip .2cm
The last part of the paper is organized in the following way. In the next section we give the proofs of 
 auxiliary results. The proofs of main results  will appear in the last section.

\section{  Proofs of auxiliary results}

The proofs of auxiliary results  will be given in the reverse order. First, we will prove Lemma 2. Then, on the basis
of Lemma 2, we will prove Lemma 1.  
\vskip .2cm
After that, we will use the C. Neumann decomposition to give the proof of Proposition~2, and then, finally we prove Proposition 1.

\vskip .2cm

\vskip .2cm
{\it Proof of Lemma 2}. We decompose the characteristic function  $\chi_{jh}$
in the operator
\begin{equation}\label{lam}
\Lambda_{2k} = 2 \sum_{j=1}^k (-1)^{j+1} a_{j,k} \, \chi_{jh}^2, \quad { a_{j,k}} = \binom{2k}{k+j}\,\binom{2k}{k}^{-1},
\end{equation}
 in the special way 
$$
\chi_{jh}  =
\frac 1j \, \sum_{i=0}^{j-1} \chi_{-(j-1)h/2+ih,h}, \quad  \chi_{t, h}(x):= \chi_h (x-t).
$$

{
In particular
\begin{align*}
\chi_{2h} =& \frac 12 \left( \chi_{-h/2,h} + \chi_{h/2,h} \right),\\
\chi_{3h} =& \frac 13 \left( \chi_{-h,h} + \chi_{0,h} + \chi_{h,h} \right),\\
\chi_{4h} =& \frac 14 \left( \chi_{-3h/2,h} + \chi_{-h/2,h} + \chi_{h/2,h} + \chi_{3h/2,h}\right), \\
\chi_{5h} =& \frac 15 \left( \chi_{-2h,h} + \chi_{-h,h} + \chi_{0,h} + \chi_{h,h} + \chi_{2h,h}\right). \\
\end{align*}

It is easy to see that
\begin{equation}\label{char}
\chi_{t, h} * \chi_{s,h} = \chi_{t+s,h}^2.
\end{equation}
Therefore,  we have by direct calculation 
\begin{equation}\label{chijh2}
\chi_{jh}^2 =\chi_{jh}*\chi_{jh} 
            = \frac 1 j \sum_{l=-j+1}^{j-1} \varphi_{l,j}  \cdot \chi_{lh,h}^2, \quad \varphi_{l,j}:= \left( 1 - \frac{|l|}j \right).
\end{equation}
\noindent
Equality \eqref{char} implies that this representation is equivalent to the following equalities for the Fej{\'er} kernel.  For $  j= 2\nu + 1, $ 
$$
\frac{1}{(2\nu+1)^2} \left( \sum_{l=-\nu}^\nu   e^{ilt} \right)^2 = \frac1{(2\nu+1)^2} \left( \frac{\sin (2\nu+1)t/2}{\sin t/2} \right)^2 = \frac 1{2\nu+1} \sum_{l=-2\nu}^{2\nu}
\varphi_{l,2\nu+1} \cdot e^{ilt},
$$
and for $ j= 2\nu$
$$
\frac{1}{(2\nu)^2} \left( \sum_{l=1}^{2\nu-1} e^{ilt/2} + e^{-ilt/2} \right)^2 = \frac1{(2\nu)^2} \left( \frac{\sin \nu t}{\sin t/2} \right)^2 = \frac 1{2\nu} 
\sum_{l=-2\nu+1}^{2\nu-1}
\varphi_{l,2\nu} \cdot e^{ilt}.
$$
}
{
\noindent
The substitution of 
$$
\chi_{jh}^2 = \frac 1 j \sum_{l=-j+1}^{j-1} \varphi_{l,j} \cdot \chi_{lh,h}^2 = 
 \frac 1{j^2} \left( j \chi_{0,h}^2 + \sum_{l=1}^{j-1} i \cdot (\chi_{(l-j)h,h}^2 + \chi_{(j-l)h,h}^2 ) \right)
 $$
into
 \eqref{lam} gives
$$
\Lambda_{2k} =  
 2\, \sum_{j=1}^{k}  \, (-1)^{j+1} a_{j,k}  \frac 1{j^2} \left( j \chi_{0,h}^2 + \sum_{l=1}^{j-1} l \cdot (\chi_{(l-j)h,h}^2 + \chi_{(j-l)h,h}^2 ) \right)
$$
\begin{equation*}
= 2 \sum_{j=1}^k (-1)^{j+1} \, a_{j,k} \, \frac 1{j} \, \chi_{0,h}^2
+ 2 \sum_{l=1}^{k-1} \chi_{lh,h}^2 \, \sum_{j=l}^k (-1)^{j+1}  \, a_{j,k}  \, \frac {j-l}{j^2}  + 
2 \sum_{l=1}^{k-1} \chi_{-lh,h}^2 \, \sum_{j=l}^k (-1)^{j+1}  \, a_{j,k}  \, \frac {j-l}{j^2}
\end{equation*}
$$
=\sum_{l =-k+1}^{k-1}  \alpha_{l,k} \, \chi_{lh,h}^2,
$$
with the coefficients 
\begin{equation}\label{aik}
\alpha_{l,k} = 2 \, \sum_{j=|l|+1}^k (-1)^{j+1} \, \frac {(j-|l|) \, a_{j,k}}{j^2}.
\end{equation}
Note that $
\mbox{sign} \, \alpha_{i,k} = (-1)^i
$
(see \cite{bks13}, Lemma A).

\noindent
By using \eqref{aik} and  the identity

\begin{equation*}\label{key}
\sigma_j:=(-1)^{j+1} j + 2 \, \sum_{l=1}^{j-1} (-1)^{l+1} \, l = \frac{1 + (-1)^{j+1}}2,
\end{equation*}
\noindent
we obtain
$$
\sum_{j=-k+1}^{k-1} |\alpha_{j,k}| = \sum_{j=-k+1}^{k-1} (-1)^j \alpha_{j,k}  =\alpha_{0,k} + 2 \sum_{l=1}^{k-1} (-1)^l \, \alpha_{l,k}
$$
$$
= 2 \, \sum_{l=1}^k (-1)^{l+1} \frac{a_{l,k}}l + 4 \, \sum_{l=1}^{k-1} (-1)^l \ 
\sum_{j=l+1}^k (-1)^{j+1} \ \frac{(j-l) \, a_{j,k}}{j^2}
$$
$$
= 2 \, \sum_{j=1}^k (-1)^{j+1} \frac{a_{j,k}}j + 2 \ \sum_{j=1}^k (-1)^{j+1} \, \frac{a_{j,k}}{j^2} \ 2 \,  \sum_{l=1}^{j-1} (-1)^l \, (j-l) 
$$
$$
= 2 \, \sum_{j=1}^k \, \frac {a_{j,k}}{j^2} \left( (-1)^{j+1}j +  2 \sum_{l=1}^{j-1} (-1)^{l+j+1} (j-l) \right)
= 2 \sum_{j=1}^k  \, \frac {\sigma_j a_{j,k}}{j^2}  = 2 \sum_{j=1,  j \,\mbox{\tiny odd}}^k \frac {a_{j,k}}{j^2} = \gamma_k.
$$
\qed
\vskip .2cm
{\it Proof of Lemma 1}  Lemma 2 implies for $ j=1, \dots, k$
{
$$
\Lambda_{2k}^{j}  =  \sum_{l=-j(k-1)}^{j(k-1)}  \alpha_{l,k}(j)  \chi_{lh, h}^{2j}, \qquad    \qquad 
\sum_{l=-j(k-1)}^{j(k-1)}  |\alpha_{l,k}(j)| \le  \gamma_{k}^{j}.
$$
}
\vskip .2cm
\noindent
{ Indeed, we have
$$
\Lambda_{2k}^2 = \Lambda_{2k}* \Lambda_{2k} = \left( \sum_{l=-(k-1)}^{k-1} \alpha_{l,k} \, \chi_{lh,h}^2 \right)*\left( \sum_{l=-(k-1)}^{k-1} \alpha_{l,k} \, \chi_{lh,h}^2 \right)
$$
$$
 = \sum_{l=-2(k-1)}^{2(k-1)} \alpha_{l,k}(2) \, \chi_{lh,h}^4, \qquad \sum_{l=-2(k-1)}^{2(k-1)} |\alpha_{l,k}(2)| \le \left(\sum_{l=-(k-1)}^{k-1} |\alpha_{l,k}| \right)^2
=\gamma_k^2,
$$
$$
\vdots
$$
$$
\Lambda_{2k}^{k}  = \sum_{l=-k(k-1)}^{k(k-1)} \alpha_{l,k}(k) \, \chi_{lh,h}^{2k}, \qquad \sum_{l=-k(k-1)}^{k(k-1)} |\alpha_{l,k}(k)| \le \left(\sum_{l=-(k-1)}^{k-1} |\alpha_{l,k}| \right)^k
=\gamma_k^k.
$$
}
 Now,  one can apply the identities  
{
$$
(g*\chi_h)'(x) = \left(h^{-1}\int_{x-h/2}^{x+h/2} f(t) \, dt \right)' = h^{-1} (f(x+h/2) - f((x-h/2)) = h^{-1} \widehat \Delta_h^1 f(x),
$$
$$
(f*\chi_h^{2})^{(2)} = ((f*\chi_h*\chi_h )')' = h^{-1}\widehat \Delta_h^1 (f*\chi_h)'  = h^{-2} \, \widehat \Delta_h^{2} f,
$$
$$
(f*\chi_h^{4})^{(4)} = ((f*\chi_h^2*\chi_h^2 )^{(2)})^{(2)} = h^{-2} \widehat \Delta_h^2 (f*\chi_h^2)^{(2)}  = h^{-4} \, \widehat \Delta_h^{4} f,
$$
$$
\vdots
$$
$$
(f*\chi_h^{2k})^{(2k)} = h^{-2k} \, \widehat \Delta_h^{2k}f,
$$
which are true a.e. for the function $f$ continuous on $(-\infty, -1) \cup (-1,1) \cup (1, +\infty)$}
and { the inequalities}
{
 $$
 \| \widehat \Delta_h^{2j} f \|_{\mathbb R} \le 4^{j}\  \| f \|_{\mathbb R}, \quad j=1, \dots, k-1,
 $$
 }
 to end the proof.}
\qed
\vskip .2cm
{\it Proof of Proposition 2} \ { By  using the Neumann decomposition of  $ g_f  $}
$$
g_f = \sum_{j=0}^{k-1} (g_f - g_f*\Lambda_{2k})*\Lambda_{2k}^j +  g_f*\Lambda_{2k}^{k}, \qquad \Lambda_{2k}^j:= \Lambda_{2k}^{j-1}*\Lambda_{2k},
$$
\noindent
{ one can estimate} the approximation of $f$ on $\mathbb I$ by algebraic polynomials $p \in \mathcal P_{n-1}$  of
 { degree $n-1 \ge 2k(2k-1)-1$:}
\begin{equation*} \label{mek}
E_{n-1}^a(f) = E_{n-1}^a(g_f)  \le \sum_{j=0}^{k-1} E_{n-1}^a ((g_f - g_f*\Lambda_{2k})*\Lambda_{2k}^j) + E^a_{n-1} (g_f*\Lambda_{2k}^k).
\end{equation*}

\noindent
We apply Theorem A for $ (g_f - g_f*\Lambda_{2k})*\Lambda_{2k}^j$,
and {the inequality}
\begin{equation}\label{in2}
\frac{n^m}{n(n-1) \cdots (n-m+1)} < 2,\quad n \ge  m(m-1), \quad m \ge 4.
\end{equation}
The  { estimate} \eqref{in2}  follows from the inequality
$$
\ln(1-x) > - 2 \ln 2 \cdot  x, \quad x \in (0,1/2).
$$
\noindent
We have

\begin{align*} \label{mek}
E_{n-1}^a(f)
&\le\sum_{j=0}^{1} \frac{\mathcal{K}_{2j}}{n^{2j}} \left\|(g_f - g_f*\Lambda_{2k})*(\Lambda_{2k}^j)^{(2j) } \right\|_{\mathbb R}
+ 2 \, \sum_{j=2}^{k-1}  \frac{\mathcal{K}_{2j}}{n^{2j}} \left\| (g_f - g_f*\Lambda_{2k})*(\Lambda_{2k}^j)^{(2j)} \right\|_\mathbb{R} \notag \\
& + 2 \, \frac{\mathcal{K}_{2k}}{n^{2k}} \, \left\| (g_f*\Lambda_{2k}^k)^{(2k)} \right\|_{\mathbb R}.
\end{align*}
\qed

\vskip .2cm
{ \it Proof of Propostion 1} 
Suppose that $p_{2k-1}^-, p_{2k-1}^+$ are the polynomials of the best approximation of $f$ on  $\mathbb I^- = [-1, -1+ 2kh]$ and $\mathbb I^+= [1-2kh,1]$
respectively.
Define 
$$
g_f(x):= \begin{cases} p_{2k-1}^+ (x), &\quad x \in (1, +\infty), \\
f(x), &\quad x \in \mathbb I, \\
p_{2k-1}^-(x), &\quad x \in (-\infty, -1). \end{cases}
$$
 Theorem B implies

\begin{equation*}\label{w2}
\sup_{x \in \mathbb I^\pm} | f(x) - p_{2k-1}^\pm (x)|
 \le \mathit{w}_{2k} \  \omega_{2k} (f, h). 
\end{equation*}

\noindent
We will prove that
\vskip .2cm
$$
W_{2k} (g_f,h)=  \| g_f - \Lambda_{2k}*g_f \|_{\mathbb R} \le  c_k \, \mathit w_{2k} \ \omega_{2k} (f, h),
$$
with 
$$
c_1 = 1, \ c_2 < 2.18, \ c_3< 2.26, \ c_4  < 2.31, \ c_k < 3, \quad k \ge 5.
$$
\noindent
By symmetry, it is sufficient to  consider only the cases 
\vskip .2cm 
1) $ x \in [0,1-kh]$,

2) $ x \in (1-kh, 1]$,

3) $ x \in (1, +\infty)$.
\vskip .2cm
\noindent
In the case 1) we have
$$
|W_{2k} (g_f, x, \chi_h^2)|  \le \binom{2k}{k}^{-1} \, \omega_{2k}(f,h).
$$

\noindent
In the second case we have  $  x+kh \ge 1$. The identity 
 $W_{2k} (p^+_{2k-1}, x, \chi_h) \equiv 0 $ yields 
{
$$
W_{2k} (g,x, \chi_h^2) = W_{2k}(g-p_{2k-1}^+,x, \chi_h^2) + W_{2k}(p_{2k-1}^+,x,\chi_h^2) = W_{2k}(g-p_{2k-1}^+,x, \chi_h^2).
$$
}
By using the  inequaltiy  \eqref{bksb}
$$
|W_{2k} (g_f - p_{2k-1}^+,x,\chi_h^2) | \le c_k \, \sup_{t \in [x-kh, x+kh]} |g_f(t) - p_{2k-1}^+(t)|,
$$
and  Whitney's theorem  \eqref{whb}
$$
\sup_{x \in I^+} | f(x) - p_{2k-1}^+(x)| \le \mathit{w}_{2k} \, \omega_{2k} (f, h),
$$
we deduce the estimate in the second case.

The third case is  similar  to the second case, when  $x \in (1, 1+kh]$.  If  $x > 1+kh$, the  $2k$--th difference is
equal to zero. Thus the estimate
$$
W_{2k} (g_f, h) \le  d_k \ \omega_{2k}(f,h),  \qquad 0<h< (2k)^{-1},
$$
is proved.

\vskip  .2cm
The proof of the { estimate}
$$
\omega_{2k} (g_f, h) \le d_k^* \ \omega_{2k} (f,h),  \qquad 0<h< (2k)^{-1},
$$
is the same as in the case $W_{2k}$. It is sufficient to consider only $x \ge 0$  and  
 instead of the inequality 
 $$
 |W_{2k}(f-p_{2k-1}^+,x, h)| \le c_k  \sup_{x \in \mathbb{I^+}} |f(x)-p_{2k-1}^+(x)|,
 $$
  to use in the cases 2) and 3) the following { inequalities}
$$ 
\left|\Delta_h^{2k} (f(x)-p_{2k-1}^\pm (x)) \right| \le (2^{2k}-1) 
 \sup_{ x \in \mathbb{I}^\pm} |f(x) - p_{2k-1}^\pm(x)| \le (2^{2k}-1) \, \mathit{w}_{2k} \, \omega_{2k} (f, h). 
$$
\qed

\vskip .2cm

\section{Proofs of main results}

\vskip .1cm

{\it Proof of Theorem 1}  Put
$$
{ \delta_k:=\delta_k (h,n):=\frac{4 \, \gamma_k}{h^2 n^2}.}
$$
Proposition 2 and Lemma 1 imply{
$${
E_{n-1}^a (f)  \le 
\left( 2 \sum_{j=0}^{k-1} \delta_k^j \, \mathcal K_{2j} - 1 -  \delta_k \mathcal K_2  \right) 
\| W_{2k}(g_f,\cdot, h) \|_{\mathbb R} 
+ 2 \mathcal K_{2k} \, 4^{-k} \, \delta_k^k \, \| \Delta_h^{2k} g_f \|_{\mathbb R}.}
$$
}
\vskip .1cm
\noindent
Now we can { apply}  Proposition 1 and  obtain{
$$
E_{n-1}^a (f)  \le 
 d_k \left( 2 \sum_{j=0}^{k}  \delta_k ^j\, \mathcal K_{2j} - 1 -   \mathcal K_2 \, \delta_k   \right) \, \omega_{2k}(f, h).
$$
}
By choosing 
$$
h = \frac{\alpha \pi}n, 
$$
and by using the identity
$$
{\sum_{j=0}^{\infty} \mathcal K_{2j} \rho^{-2j} = \sec \left( \frac \pi{2 \rho} \right), \quad \rho :=  { \delta_k^{-1/2} = \alpha  \ \ \sqrt{\frac{\pi^2 }{4\gamma_k}}  }  > \alpha  >1,}
$$
we deduce the estimate
$$
J(2k, \alpha) { \le 3 \, (2 + e^{-2}) \, \left(2 \sec \left( \frac \pi{2 \rho} \right) - 1 -\frac {\mathcal K_2} {\rho^2} \right)} \le 3 \, (2 + e^{-2}) \, \left(2 \sec \left( \frac \pi{2 \alpha} \right) - 1 -\frac {\mathcal K_2} {\alpha^2} \right), \quad \alpha >1.
$$
{ The last estimate follows from the fact that  $ 2 \sec \left( \frac \pi{2 x}\right) - 1 -\frac {\mathcal K_2}{x^2} $ is a decreasing function for $x>1$.}
So, the upper estimate for $J(2k, \alpha)$  is proved.

\vskip .1cm

At last  we  prove the lower estimate. Consider the function 
$$
f_0(x) = \begin{cases}1, \quad & x = -1,\\ 
                      0, \quad & x \in \mathbb I, \ x \ne -1. \end{cases}
$$
The best approximation of this function is $\ge 1/2$. But the $k$--th difference 
$$
 \Delta_h^k f(x) = \sum_{j=0}^k (-1)^{k-j} \binom{k}j f(x +jh), \quad h \in (0, (1-x)/k), 
$$
for  $x= -1$ is equal to $(-1)^k$.

\noindent
Consider now a regularization of $f_0$.
Suppose  that $k \ge 2, \ \varepsilon \in (0,1/k)$ and define
$$
f_{\varepsilon,k} (x) = f_{\varepsilon} (x):=\begin{cases} \frac{(-1)^{k-1}}{\varepsilon^{k-1}} \,
(1 + x - \varepsilon)^{k-1},& \quad  -1   \le x \le -1+\varepsilon, \\
0,& \quad  -1 +\varepsilon < x \le 1. \end{cases}.
$$
We will use the representation
$$
\Delta_h^k f_{\varepsilon} (x) = \int_0^h \, du_1 \, \cdots \, \int_0^h \Delta_h^1 f_{\varepsilon} ^{(k-1)} (x + u_1 + \cdots + u_{k-1} ) \, d u_{k-1}.
$$
Note that
$$
\left| f_{\varepsilon}^{(k-1)}(x) \right| = \begin{cases} \frac{(k-1)!}{\varepsilon^{k-1}}, \qquad &x \in (-1, -1+\varepsilon),\\
0, \qquad &x \in (-1+\varepsilon, 1).
\end{cases}
$$
{ For $h>0$  the first difference}
$$
\Delta_h^1 f_{\varepsilon} ^{(k-1)} (x + u_1 + \cdots + u_{k-1} ) = f_{\varepsilon} ^{(k-1)} (x + h+ u_1 + \cdots + u_{k-1} )-
f_{\varepsilon} ^{(k-1)} (x + u_1 + \cdots + u_{k-1} )
$$
is not equal to zero only if 
$$
0 \le u_1 + \cdots + u_{k-1} \le \varepsilon.
$$
\noindent
From
$$
\left| \Delta_h^1 f^{(k-1)}(x) \right| \le \frac{(k-1)!}{\varepsilon^{k-1}},
$$
\noindent
and 

$${\int \cdots \int}_{0 \le u_1 + \cdots + u_{k-1} \le \varepsilon} \, d u_1 \cdots d u_{k-1} = \frac{\varepsilon^{k-1}}{(k-1)!},
$$
we obtain for 
 $x, x + kh \in  \mathbb I $ 

$$
\left| \Delta_h^k f_\varepsilon (x) \right| \le 1.
$$
It is clear, that for small $\varepsilon >0$
$$
E_{n-1}^a (f_\varepsilon)  \ge 1/2 - \delta(\varepsilon), \quad  \lim_{\varepsilon \to 0} \delta(\varepsilon) = 0. $$
{ In the  case $k=1$  sufficient to consider the function $f_{\varepsilon,2}$}.
\qed

{\it Proof of Theorem 2} \ \ In the  case $k=1$ we have
{
$$
E_{n-1}^a(f) \le W_2(g_f, h) + E_{n-1}^a (g_f*\chi_h^2) \le  0.5 \, \omega_2 (g_f,h) + \mathcal K_2 \, n^{-2} \| (g_f*\chi^2_h)^{(2)} \| \le (0.75 + 1.5\, \mathcal K_2 (hn)^{-2} ) 
\, \omega_2 (f, h).
$$
 In this case  the inequalities (see Proposition 1)
$$
W_2(g_f,h) \le 0.5  \ \omega(g_f,h) \le 3/4 \ \omega_2 (f, h)
$$
are better than inequality 
$$
W_2(g_f,h) \le \omega_2(f,h).
$$

In the  case $ k \ge 2$ the estimates $W_{2k}(g_f,h) \le  d_k \ \omega_{2k} (f,h)$   give  better results.}

If $k=2,3,4,$ then  \ Theorem 2 follows from { the inequality  (see the proof of Theorem 1) 
$$
E_{n-1}^a (f)  \le 
 \left( 2  d_k \, \sum_{j=0}^{k-1} \left(  \delta_k^j\, \mathcal K_{2j} - 1 -   \delta_k \mathcal K_2  \right)  + 2\, d_k^* \,  \mathcal K_{2k} \, 4^{-k} \, \delta_k^k  \right) \, \omega_{2k}(f, h),
$$
}
\noindent
and { the estimates} of $d_k, d_k^*, \gamma_k $   (see Proposition 1 and definition of $\gamma_k$). 
~\qed

\noindent
\begin{minipage}[t]{6.5cm}
\scriptsize{
Alexander Babenko \\
Institute of Mathematics and Mechanics \\
Ural Branch of the Russian\\ Academy of Sciences \\
S.\,Kovalevskoi Str. 16, \\ Ekaterinburg, 620990, Russia,\\
{Ural Federal University,} \\
{Ekaterinburg, Russia}\\
\textrm{babenko@imm.uran.ru}}
\end{minipage}
\begin{minipage}[t]{4.5cm}
\scriptsize{Yuriy Kryakin \\ Institute of Mathematics \\
University of Wroclaw \\ Plac Grunwa{l}dzki 2/4, \\ 50-384
Wroclaw, Poland\\ \textrm{kryakin@gmail.com,
kryakin@math.uni.wroc.pl}}
\end{minipage}

\end{document}